\documentstyle[12pt]{article}
\begin{document}
\title{
Quantum Dynamical $R$-matrices and Quantum Frobenius Group
}
\author{G.E.Arutyunov
\thanks{Steklov Mathematical Institute,
Vavilov 42, GSP-1, 117966, Moscow, Russia; arut@class.mi.ras.ru}\\
and \\
S.A.Frolov
\thanks{Steklov Mathematical Institute,
Vavilov 42, GSP-1, 117966, Moscow, Russia; frolov@class.mi.ras.ru}\\
}
\date {}
\maketitle
\begin{abstract}
We propose an algebraic scheme for quantizing the
rational Ruijsenaars-Schneider model in the $R$-matrix formalism.
We introduce a special parameterization of the cotangent bundle over
$GL(N,{\bf C})$. In new variables the standard symplectic structure
is described by a classical (Frobenius) $r$-matrix and by a new
dynamical $\bar{r}$-matrix. Quantizing both of them we find the quantum 
$L$-operator algebra and construct its particular representation 
corresponding to the rational Ruijsenaars-Schneider system. 
Using the dual parameterization of the cotangent bundle we also derive 
the algebra for the $L$-operator of the trigonometric Calogero-Moser 
system.   
\end{abstract} 
\newpage 
\section{Introduction} 
As soon as the classical dynamical $r$-matrices first appeared \cite{AT} on 
the scene of integrable many body systems, the problem of their 
quantization became of a real interest. The main hope related to 
this problem is to find a new algebraic structure that ensures the 
integrability of the corresponding quantum models.

We recall \cite{BV} that having a finite-dimensional completely 
integrable system with the Lax representation $\frac{dL}{dt}=[M,L]$
one can always write the Poisson algebra of $L$-operators in the 
$r$-matrix form. However, in general, an $r$-matrix appears to be 
a nontrivial function of dynamical variables.
At present the classical dynamical $r$-matrices are known for the
rational, trigonometric \cite{AT,ABT} and elliptic \cite{Skl,BS} 
Calogero-Moser (CM) systems, as well as for their relativistic 
generalizations -- rational, trigonometric \cite{AR,Sur} and elliptic 
\cite{NKSR,Surel} Ruijsenaars-Schneider (RS) systems \cite{R}. 

The problem of quantizing the dynamical $r$-matrices is rather nontrivial
since, in general, such $r$-matrices do not satisfy a single
closed  equation of the Yang-Baxter type, 
from which they can be uniquely determined. 
Up to now there exists only one example of a quantum dynamical 
$R$-matrix related to the quantum spin CM system \cite{ABB}.  
This $R$-matrix solves the Gervais-Neveu-Felder equation 
\cite{GN,Fel} and has a nice 
interpretation in terms of quasi-Hopf algebras \cite{BBB}.

A natural way to understand the origin of dynamical $r$-matrices
is to consider the Hamiltonian reduction procedure \cite{ABT,AM}.
Factorizing a free motion on an initial phase space by
the action of some symmetry group we get nontrivial dynamics on
the reduced space. The $r$-matrix appears in the Dirac bracket, which 
describes the phase structure of the reduced space.
In our recent papers \cite{AFM,AFM1} we obtained 
the elliptic RS model, being
the most general one among the integrable systems of the CM
and RS types, by using two different reduction schemes.
In the first scheme, the affine Heisenberg double was used
as the initial phase space and in the second one we considered
the cotangent bundle over the centrally extended group of double loops.

The aim of this paper is to quantize the reduction scheme leading
to the dynamical systems of the RS type. Although the
most interesting is the spectral-dependent elliptic case \cite{AFM,AFM1},
to clarify the general approach in this paper we restrict ourselves 
to considering the simplest rational model. 

Our construction is based on a special parameterization of the cotangent 
bundle over the group $GL(N,{\bf C})$. This parameterization is similar to 
the one considered in \cite{AF}. However, instead of Euler angles, we use 
another system of generators to parameterize the "momentum". These
generators obey a quadratic Poisson algebra describing by 
two $r$-matrices $r$ and $\bar{r}$. The matrix $r$ solves the 
$N$-parametric classical Yang-Baxter equation and is related to the 
special Frobenius subgroup in $GL(N,{\bf C})$.

We define a special matrix function $L$ on $T^*G$ invariant
with respect to the action of this Frobenius subgroup. We call
this function the "$L$-operator" since the Poisson algebra
of $L$ literally coincides with the one for the 
rational RS model \cite{Sur}. Moreover, performing the 
Hamiltonian reduction leading to the RS model \cite{GN}, 
one can recognize in $L$ the standard $L$-operator of the rational RS model.

Then we pass to the quantization. The quadratic Poisson algebra can 
be quantized by using the $R$-matrix approach \cite{F,KS}. The 
compatibility of the corresponding quantum algebra implies the 
quantum Yang-Baxter equation for $R$ and some new equations involving 
$R$ and $\bar{R}$. We solve these equations and get an explicit
form for $R$ and $\bar{R}$. Coming back to the original generators 
of $T^*G$ we recover the standard commutation relations of the 
quantum cotangent bundle.

We derive a new quadratic algebra which is satisfied by the "quantum"
$L$-operator. The matrices $R$ and $\bar{R}$ come in this algebra
in a nontrivial way. This explains the nature of the dynamical
$r$-matrices (classical and quantum) as the composite objects
constructed from the more elementary blocks $R$ and $\bar{R}$. 

It follows from our construction that the quantum $L$-operator is 
factorized in the form $L=WP$. Here $W$ satisfies
the defining relations of the quantum Frobenius group, 
$W_1W_2R_{12}=R_{12}W_2W_1$, where $R$, being the quantization of $r$, 
is the $N$-parametric solution to the quantum Yang-Baxter equation.
The diagonal matrix $P$ plays the role of a generalized
momentum. We find the simplest representation of the $L$-operator
algebra and relate it with the rational RS model.

The trigonometric CM system is known to be dual
to the rational RS model \cite{R}. This duality is explained
by the existence of the dual parameterization of $T^*G$. In this
parameterization the CM model can be easily quantized and we
get the commutation relations satisfied by the corresponding quantum 
$L$-operator. 

\section{Frobenius algebra and dynamical $r$-matrices}
In this section we introduce a special parameterization of the cotangent
bundle $T^*G$ over the matrix group $G=GL(N,{\bf C})$. As a manifold the 
space $T^*G$ is naturally isomorphic to ${\cal G}^*\times G$, where
${\cal G}^*$ is dual to the Lie algebra  
${\cal G}=\mbox{Mat}(N,{\bf C})$. The standard Poisson structure on 
$T^*G$ can be written in terms of variables $(A,g)$, where $A\in 
{\cal G}^{*}\approx {\cal G}$ and $g\in G$, as follows 
\begin{eqnarray} 
\{ A_1,A_2\} &=&\frac 1 2 [C,A_1 
-A_2] \label{AA}\\
\{ A_1,g_2\} 
&=&g_2 C \label{pb} \\
\{g_1,g_2\} 
&=& 0\label{gg}
\end{eqnarray} 
Here we use the standard tensor notation 
and $C=\sum_{i,j}E_{ij}\otimes E_{ji}$ is the permutation operator.  

Any matrix $A$ belonging to an orbit  of maximal 
dimension in ${\cal G}^*$ admits the factorization: 
\begin{equation}
A=TQT^{-1},
\label{facA}
\end{equation}
where $Q$ is a diagonal matrix with entries $q_i,~q_i\neq q_j$.
\noindent
We also fix the order of $q_i$: $q_1<q_2 \ldots 
<q_n$ by using the action of the Weyl group.  It is obvious that the 
matrix $T$ in eq.(\ref{facA}) is not uniquely defined. Indeed, one can 
multiply $T$ by an arbitrary diagonal matrix from the right. We remove 
this ambiguity by imposing the following condition: 
\begin{equation} Te=e, 
\label{Te}
\end{equation}
where $e$ is a vector with all $e_i=1$. 
\noindent
The choice of (\ref{Te}) is motivated by the study
of the reduction procedure leading to the Calogero-type integrable 
systems \cite{OP1}.
Let us note that the condition (\ref{Te}) defines a Lie subgroup 
$F\subset G$. The corresponding Lie algebra $\cal F$ has the natural basis
$F_{ij}=E_{ii}-E_{ij}$, where $E_{ij}$ are the standard matrix unities.
The commutation relations of $F_{ij}$ are
$$
[F_{ij},F_{kl}]=\delta_{ik}(F_{il}-F_{ij})+
\delta_{il}(F_{kj}-F_{kl})+\delta_{jk}(F_{ij}-F_{il}).
$$
It is worthwhile to mention that $\cal F$ is not only the Lie algebra but 
also an associative algebra with respect to the usual matrix multiplication:
$$
F_{ij}F_{kl}=\delta_{ik}F_{il}+
\delta_{jk}(F_{ik}-F_{il}).
$$

Let us now rewrite the Poisson structure (\ref{AA}) in terms of 
the variables
$T$ and $Q$. It is obvious that the coordinates $q_i$ commute with $A$,
$T$ and $Q$ since they belong to the center of the Poisson algebra 
(\ref{AA}).  Thus, it is enough to calculate the Poisson bracket $\{T,T\}$. 
We have 
$$ \{T_{ij},T_{kl}\}=\sum_{mn,ps}\frac{\delta 
T_{ij}}{\delta A_{mn}} \frac{\delta T_{kl}}{\delta 
A_{ps}}\{A_{mn},A_{ps}\}.  $$ 
To find $\frac{\delta T_{ij}}{\delta 
A_{mn}}$ we perform the variation of (\ref{facA}):  $$ T^{-1}\delta T Q-Q 
T^{-1}\delta T + \delta Q = T^{-1}\delta A T.  $$ 
This equation can be 
easily solved, and we obtain the derivatives 
\begin{equation} 
\frac{\delta 
T_{ij}}{\delta A_{mn}}= \sum_{a\neq j}\frac{1}{q_{ja}} 
(T_{ia}T_{nj}T_{am}^{-1}+T_{ij}T_{na}T_{jm}^{-1})
\label{dTdA}
\end{equation}
and                         
\begin{equation}
\frac{\delta q_{i}}{\delta A_{mn}}=T_{ni}T_{im}^{-1},
\label{dQdA}
\end{equation}
where $q_{ij}\equiv q_i-q_j$.

By using (\ref{dTdA}) we get
\begin{equation}
\{T_1,T_2\}=T_1T_2 r_{12}(q),
\label{TT}
\end{equation}
where the $r$-matrix
\begin{equation}
r_{12}(q) = \sum_{i\neq j}\frac{1}{q_{ij}}F_{ij}\otimes F_{ji}
\label{r}
\end{equation}
appears. 
It is clear, that $r_{12}(q)$ should be a skew-symmetric 
solution of the classical Yang-Baxter equation (CYBE).
The origin of this $r$-matrix can be easily understood if we note
that $\cal F$ is a Frobenius Lie algebra,
i.e. there is a nondegenerate 2-cocycle (coboundary)
on ${\cal F}$:
\begin{equation}
\omega(X,Y)=\mbox{tr}(Q[X,Y]), ~~~X,Y \in {\cal F}.
\label{coc}
\end{equation}
According to \cite{BD}, to any Frobenius Lie algebra one can associate 
a skew-symmetric solution of the CYBE by inverting the corresponding
2-cocycle. One can check that the cocycle (\ref{coc}) corresponds to
$r_{12}(q)\in {\cal F}\wedge {\cal F}$. Coming back to (\ref{facA})
we see that any orbit of maximal dimension 
in ${\cal G}^*$ can be supplied with the structure of the Frobenius group. 
It is worthwhile to note that $\omega$ is the Kirillov symplectic 
form on the coadjoint orbit of the maximal dimension parameterized
by $Q$.
 
Now, following \cite{AF}, we introduce a special parameterization for the group
element $g$. Namely, we take $A'=gAg^{-1}$, which Poisson commutes
with $A$ and possesses the Poisson bracket
$$
\{ A_1',A_2'\} =-\frac 1 2 [C,A_1' -A_2']
$$

Diagonalizing it with the help of the matrix $U$, $Ue=e$,
$A'=UQU^{-1}$, we find that 
\begin{equation}
\label{g}
g=UPT^{-1},
\end{equation}
where $P$ is some diagonal matrix. 
\noindent It is obvious that the Poisson bracket for $U$ is given by
\begin{equation}
\{U_1,U_2\}=-U_1U_2r_{12}(q)
\label{UU}
\end{equation}
and that $\{T,U\}=0$. To proceed with the calculation of the
brackets $\{U,P\}$, $\{T,P\}$, $\{P,P\}$ and $\{P,Q\}$ we should
use the derivative $\frac{\delta U_{ij}}{\delta A_{mn}'}$
which is given by (\ref{dTdA}) with  the replacement $T\to U$, $A\to A'$.
After simple calculations we get
\begin{eqnarray}
\{P_1,P_2\}&=&0, \nonumber\\
\{Q_1,P_2\}&=&P_2 \sum_i E_{ii}\otimes E_{ii}.
\label{QP}
\end{eqnarray}
Introducing $p_i=\log P_i$ we conclude that $\{q_i,p_j\}=\delta_{ij}$.
We also find the remaining Poisson brackets
\begin{eqnarray}
\label{UP}
\{U_1,P_2\}&=&U_1P_2\bar{r}_{12}(q),\\
\{T_1,P_2\}&=&T_1P_2\bar{r}_{12}(q),
\label{TP}
\end{eqnarray}
where we have introduce a new matrix
\begin{equation}
\bar{r}_{12}(q)=\sum_{i\neq j}\frac{1}{q_{ij}}F_{ij}\otimes E_{jj}.
\label{rbar}
\end{equation}
The Jacobi identity leads to a set of equations on the matrices $r$
and $\bar{r}$. However, we postpone the discussion of these equations 
till the next section, where the quantization of $T^*G$ will be 
given in terms of variables $Q,~T,~P,~U$.

Let us define the $L$-operator as the following function of phase
variables, being invariant under the action of the Frobenius group:
\begin{equation}
L=T^{-1}gT = T^{-1}UP.
\label{L}
\end{equation}
By using (\ref{TT}), (\ref{UU}) and (\ref{QP}-\ref{TP})
one can easily find the Poisson brackets containing $L$:
\begin{eqnarray}
\label{QL}
\{Q_1,L_2\}&=&L_2\sum_i E_{ii} \otimes E_{ii}, \\
\label{TL}
\{T_1,L_2\}&=&T_1L_2\bar{r}_{12}(q)-T_1r_{12}(q)L_2, \\
\label{LL}
\{L_1,L_2\}&=&r_{12}(q)L_1L_2+
L_1L_2(\bar{r}_{12}(q)-\bar{r}_{21}(q)-r_{12}(q))\\
\nonumber
&+& L_1\bar{r}_{21}(q)L_2-L_2\bar{r}_{12}(q)L_1.
\end{eqnarray}

We see that the brackets for $L$ and $Q$ are the brackets 
for the $L$-operator and coordinates of the rational RS model
found in \cite{Sur}. It is obvious that  
$I_n=\mbox{tr}L^n=\mbox{tr}g^n$ form a set of mutually commuting 
functions. 

Let us note that the $L$-operator (\ref{L}) has the form $L=WP$,
where $W=T^{-1}U$. Since both $T$ and $U$ are elements of the 
Frobenius group $F$, $W$ also belongs to $F$. Calculating
the Poisson bracket for $W$ we see that it coincides with
the Sklyanin bracket defining on $F$ the structure of a Poisson-Lie
group:
\begin{equation}
\{W_1,W_2\}=[r_{12}(q),W_1W_2].
\label{WW}
\end{equation}
The Poisson relations of $W$ and $P$ are
\begin{equation}
\{W_1,P_2\}=-P_2[\bar{r}_{12}(q),W_1].
\label{WP}
\end{equation}

The well-known property of the Poisson algebra (\ref{WW})
is the existence of a family of mutually commuting functions
$J_n=\mbox{tr}W^n$. Moreover, it turns out that $J_n$ commute
not only with themselves but also with $P$ and $Q$. In Section 3
we show that the same property holds in the quantum case.

Now we construct the simplest representation of the Poisson algebra
(\ref{WW}), (\ref{WP}) and (\ref{QP}) and relate it with the rational
RS model. In fact, this representation corresponds to the zero
dimensional symplectic leaf of the bracket (\ref{WW}).

To this end we employ the Hamiltonian reduction procedure.
We recall that $G$ acts on $T^{*}G$ by the transformations
$A\to hAh^{-1}$, $g\to hgh^{-1}$ in a Hamiltonian way. The corresponding
moment map $\mu$ has the form $\mu=gAg^{-1}-A$. 
Performing the Hamiltonian
reduction we fix it's value to be
$$
gAg^{-1}-A=-\gamma (e\otimes e - 1),
$$ 
where $\gamma$ is an arbitrary constant.
In terms of $(T,L,Q)$ variables this equation acquires the form
\begin{equation}
T(LQL^{-1}-Q)T^{-1}=-\gamma (e\otimes e - 1).
\label{mom}
\end{equation} 
Since $Te=e$ and $L=WP$ the last equation can be written as 
\begin{equation}
WQ-QW-\gamma W=-\gamma e\otimes(eU).
\label{mom1}
\end{equation} 
Eq.(\ref{mom1}) can be elementary solved and one gets
\begin{equation}
W=\sum_{i,j}\frac{\gamma}{\gamma+q_{ij}} e_ib_jE_{ij},
\label{mom2}
\end{equation} 
where $b=eU$. If we recall that $W$ should be an element of $F$, i.e.
$We=e$, then we find the coefficients $b_j$:
\begin{equation}
b_j=\frac{1}{\gamma}\frac{\prod_{a}(q_{aj}+\gamma)}{\prod_{a\neq j} q_{aj}}
\label{b}
\end{equation} 
and thereby
\begin{equation}
W_{ij}=\frac{\prod_{a\neq i}(q_{aj}+\gamma)}{\prod_{a\neq j}q_{aj}}. 
\label{W} 
\end{equation} 

One can check that $W$ given by (\ref{W}) has the desired Poisson 
brackets (\ref{WW}) and (\ref{WP}). We do not give an explicit
proof of this statement since in section 3 we show that the same 
function $W$ realizes the representation for the corresponding 
quantum algebra.  

The relation of our $L$-operator with the standard Ruijsenaars 
$L$-operator \cite{R} is given by the following canonical 
transformation:
$$
q_i\to q_i,~~~P_i\to 
\prod_{a\neq i}\frac{(q_{ai}-\gamma)^{1/2}}{(q_{ai}+\gamma)^{1/2}}P_i. 
$$

\section{Quantum $L$-operator algebra}
In this section we shall quantize $T^{*}G$ in terms of variables
$Q,~T,~P,~U$ and obtain the permutation relations for the 
quantum $L$-operator. The algebra of functions on $T^*G$
can be easily quantized and one gets an associative algebra 
generated by $A$ and $g$ subject to the standard relations:

\begin{eqnarray} 
\left[A_1,A_2\right]&=&\frac{1}{2} {\hbar} [C,A_1 - A_2] \label{qAA}\\
\left[A_1,g_2\right]&=&{\hbar} g_2 C \label{qAg} \\     
\label{qgg} 
\left[g_1,g_2\right]&=& 0, 
\end{eqnarray} 
where $\hbar$ is a quantization parameter.
The commutation relations for the generators $T,Q,U,P$ can be 
straightforwardly written by using the ideology of the Quantum 
Inverse Scattering Method \cite{F,KS} (we present only nontrivial relations):  
\begin{eqnarray} 
\label{qTT} 
T_1T_2&=&T_2T_1R_{12}(q),~~U_1U_2=U_2U_1R_{12}^{-1}(q)\\
\label{qTP} 
T_1P_2&=&P_2T_1\bar{R}_{12}(q),~~U_1P_2=P_2U_1\bar{R}_{12}(q)\\
\label{qQP}
&&[Q_1,P_2]=\hbar P_2 \sum_i E_{ii}\otimes E_{ii}.
\end{eqnarray}
Here $R(q)$ and $\bar{R}(q)$ are quantum dynamical $R$-matrices
having the following behavior near $\hbar =0$:
\begin{eqnarray}
\nonumber
R(q)=1+\hbar r(q) + o(\hbar),~~\bar{R}(q)=1+\hbar \bar{r}(q) + o(\hbar).
\end{eqnarray}
It follows from the compatibility conditions that the $R$-matrices 
should satisfy the following set of equations
\begin{eqnarray}
\label{RR}
R_{12}(q)R_{21}(q)&=&1,\\
\label{RRR}
R_{12}(q)R_{13}(q)R_{23}(q)&=&R_{23}(q)R_{13}(q)R_{12}(q),\\ 
\label{RRbRb}
R_{12}(q)\bar{R}_{13}(q)\bar{R}_{23}(q)&=&
\bar{R}_{23}(q)\bar{R}_{13}(q)P_3^{-1}R_{12}(q)P_3,\\ 
\label{RbRb}
\bar{R}_{12}(q)P_2^{-1}\bar{R}_{13}(q)P_2&=&
\bar{R}_{13}(q)P_3^{-1}\bar{R}_{12}(q)P_3
\end{eqnarray}
Let us demonstrate how one gets, for example, eq.(\ref{RRbRb}). 
This equation follows from (\ref{qTP}) and the following chain of relations:
$$
T_1T_2P_3=T_1P_3T_2\bar{R}_{23}=P_3T_1\bar{R}_{13}T_2\bar{R}_{23}=
P_3T_2T_1R_{12}\bar{R}_{13}\bar{R}_{23}=
T_2P_3\bar{R}^{-1}_{23}T_1R_{12}\bar{R}_{13}\bar{R}=
$$
$$
T_2T_1P_3\bar{R}^{-1}_{13}\bar{R}^{-1}_{23}R_{12}\bar{R}_{13}\bar{R}=
T_1T_2R^{-1}_{12}P_3\bar{R}^{-1}_{13}\bar{R}^{-1}_{23}R_{12}\bar{R}_{13}\bar{R}.
$$

The solution of the Yang-Baxter equation (\ref{RRR}) can be easily found 
if one notes that $r^2_{12}=0$. It is known 
(see e.g. \cite{Chari} Prop.6.4.13) that if
$r\in \mbox{Mat}(N,{\bf C})\otimes \mbox{Mat}(N,{\bf C})$ satisfies
$r^3=0$ and solves the classical Yang-Baxter equation then $R=e^r$ is
a solution of QYBE. Therefore, 
$R_{12}=e^{\hbar r_{12}}=1+\hbar r_{12}$ is a desired solution of the 
QYBE. 

The solution of eq.(\ref{RRbRb}) can be found if one supposes that
$\bar{R}$ has the same matrix structure  as $\bar{r}$ does:
\begin{equation}
\label{RbA}
\bar{R}_{12}(q)=1+\hbar\sum_{i\neq 
j}\bar{r}_{ij}(\hbar,q)F_{ij}\otimes E_{jj} 
\end{equation} Then the 
following $\bar{R}$-matrix is a solution of eqs.(\ref{RRbRb}) and 
(\ref{RbRb}):  
\begin{equation} \label{Rb} 
\bar{R}_{12}(q)=1+
\sum_{i\neq j}\frac{\hbar}{q_{ij}-\hbar}F_{ij}\otimes E_{jj}.
\end{equation}
It is not difficult to verify that 
$
\bar{R}_{12}^{-1}(q)=1-
\sum_{i\neq j}\frac{\hbar}{q_{ij}}F_{ij}\otimes E_{jj}.
$

Now we should show that the generators $A=TQT^{-1}$ and $g=UPT^{-1}$
satisfy the commutation relations (\ref{qAA}-\ref{qgg}). From
the relations (\ref{qTT}-\ref{qQP}) we get
\begin{eqnarray}
\nonumber
\left[A_1,A_2\right]&=&T_2T_1 
(R_{12}Q_1R_{21}Q_2-Q_2R_{12}Q_1R_{21})T_1^{-1}T_2^{-1}\\ \nonumber 
\left[A_1,g_2\right]&=&g_2T_2T_1 
(R_{12}(Q_1+\hbar \sum_i E_{ii}\otimes E_{ii})\bar{R}^{-1}_{12}R_{12}-R_{12}Q_1)T_2^{-1}T_1^{-1}  \\
\nonumber
\left[g_1,g_2\right]&=&U_2U_1 
(R^{-1}_{12}P_1\bar{R}^{-1}_{21}P_2\bar{R}^{-1}_{12}R_{12}-
P_2\bar{R}^{-1}_{12}P_1\bar{R}^{-1}_{21})T_2^{-1}T_1^{-1}
\end{eqnarray}
By using the following identities that can be checked by direct 
computation
\begin{eqnarray}
\nonumber
R_{12}Q_1R_{21}Q_2-Q_2R_{12}Q_1R_{21}&=&\hbar (CQ_1R_{21}-R_{12}Q_1C)\\ 
\nonumber
R_{12}(Q_1+\hbar \sum_i E_{ii}\otimes E_{ii})\bar{R}^{-1}_{12}R_{12}
-R_{12}Q_1 &=&\hbar C, \\
\label{RRbPPRb}
R^{-1}_{12}P_1\bar{R}^{-1}_{21}P_2\bar{R}^{-1}_{12}R_{12}-
P_2\bar{R}^{-1}_{12}P_1\bar{R}^{-1}_{21}&=&0
\end{eqnarray}
one derives the desired commutation relations for $A$ and $g$.
Let us finally present the permutation relations for $Q$, $T$ and $L$,
which can be easily obtained by using (\ref{qTT}-\ref{qQP}) and 
(\ref{RRbPPRb}):  
\begin{eqnarray}
\label{QTL}
[Q_1,L_2]&=&\hbar L_2 \sum_i E_{ii}\otimes E_{ii}, \\
\label{tl}
T_1R_{12}L_2&=&L_2T_1\bar{R}_{12}\\
\label{qLL}
L_1\bar{R}^{-1}_{21}L_2\bar{R}^{-1}_{12}R_{12}\bar{R}_{21}&=&
R_{12}L_2\bar{R}^{-1}_{12}L_1
\end{eqnarray}

It is clear that just as in the classical case
the quantities $I_n=\mbox{tr}g^n$
form a set of mutually commuting operators.  
Let us show that $I_n$ can be expressed in terms of $L$ and $Q$
solely and thereby they can be interpreted as the 
quantum integrals of motion for the rational RS model.
By using the definition of $L$ we rewrite $I_n$ as
$$
I_n=\mbox{tr}g^n=\mbox{tr}TL^nT^{-1}=\mbox{tr}_{12}C_{12}T_1L_1^nT_2^{-1}=
\mbox{tr}_{12}C_{12}^{t_2}T_1L_1^n\stackrel{t~~}{T_2^{-1}}
$$
where $t_2$ denotes the matrix transposition in the second factor of 
the tensor product.  

It follows from (\ref{tl}) that 
$$ 
L_1\stackrel{t~~}{T_2^{-1}}=\stackrel{t~~}{T_2^{-1}}
R_{12}^{t_2}L_1\bar{R}_{21}^{t_2}.
$$
Applying this relation we derive
$$
I_n=\mbox{tr}_{12}C_{12}^{t_2}T_1\stackrel{t~~}{T_2^{-1}}
R_{12}^{t_2}L_1\bar{R}_{21}^{t_2} \cdots R_{12}^{t_2}L_1\bar{R}_{21}^{t_2}.
$$
Exchanging $T_1$ and $\stackrel{t~~}{T_2^{-1}}$ with the help of 
eq.(\ref{qTT}) one gets
$$
I_n=\mbox{tr}_{12}C_{12}^{t_2}\stackrel{t~~}{T_2^{-1}}T_1
L_1\bar{R}_{21}^{t_2} \cdots R_{12}^{t_2}L_1\bar{R}_{21}^{t_2}. 
$$
Since $C_{12}^{t_2}\stackrel{t~~}{T_2^{-1}}T_1=C_{12}^{t_2}$
and $\bar{R}_{21}^{t_2}C_{12}^{t_2}=C_{12}^{t_2}$ we finally arrive at
$$
I_n=
\mbox{tr}_{12}C_{12}^{t_2}
L_1\bar{R}_{21}^{t_2}R_{12}^{t_2}
L_1\bar{R}_{21}^{t_2}R_{12}^{t_2}L_1\cdots
L_1\bar{R}_{21}^{t_2}R_{12}^{t_2}L_1. 
$$
It is natural to regard this expression for $I_n$ as a "quantum trace" 
of the operator $L^n$.

Just as in the classical case the quantum $L$-operator has the 
form $L=WP$, where $W=T^{-1}U$ satisfies the defining relations
of the quantum Frobenius group:
\begin{equation}
R_{12}W_2W_1=W_1W_2R_{12}.
\label{qWW}
\end{equation}
The algebra (\ref{QTL}), (\ref{qLL}) rewritten in terms of $Q$, $P$ and 
$W$ is given by (\ref{qQP}), (\ref{qWW}) and by the relation
\begin{equation}
W_1 P_2 \bar{R}_{12}^{-1}=P_2 \bar{R}_{12}^{-1} W_1.
\label{qWP}
\end{equation}
This shows that the representation theory for $L$ essentially
reduces to the one for the quantum Frobenius group.

It is known that the algebra (\ref{qWW}) admits
a family of mutually commuting operators given by \cite{M}:
$$
J_n=\mbox{tr}_{1\ldots n}\left[
\hat{R}_{12}\hat{R}_{23}\ldots 
\hat{R}_{n-1,n}W_1\ldots W_n \right],
$$
where $\hat{R}_{ij}=R_{ij}C_{ij}$. 
\noindent
Now we demonstrate that $J_n$ commutes with $P$. For the sake
of clarity we do it for $n=3$. It follows from (\ref{qWP}) that
$$
J_3P_4=
\hat{R}_{12}\hat{R}_{23}P_4 
\bar{R}_{14}^{-1}W_1\bar{R}_{14}
\bar{R}_{24}^{-1}W_2\bar{R}_{24}
\bar{R}_{34}^{-1}W_3\bar{R}_{34}.
$$
Equation (\ref{RRbRb}) written in terms of $\hat{R}$ acquires the form
$$
P_3^{-1}\hat{R}_{12}(q)P_3=
\bar{R}_{23}^{-1}\bar{R}_{13}^{-1}
\hat{R}_{12}\bar{R}_{13}\bar{R}_{23}.
$$
By using this equation we can push $P_4$ on the left. Taking into
account that $\bar{R}_{12}$ is diagonal in the second space we get
$$
J_3P_4=P_4\bar{R}_{14}^{-1}\bar{R}_{24}^{-1}\bar{R}_{34}^{-1}
\hat{R}_{12}\hat{R}_{23}W_1W_2W_3\bar{R}_{34}\bar{R}_{24}\bar{R}_{14}.
$$
Taking the trace in the first, second and third spaces one gets
the desired property.

Now we give an example of the representation of the algebra
(\ref{qWW}) and (\ref{qWP}). Namely, we prove that the $W$-operator given 
by eq.(\ref{W}) realizes this algebra with 
$P_j=e^{-\hbar \frac{\partial}{\partial q_j}}$.

It is obvious that $[W_1,W_2]$ should be equal to zero 
since $W$ depends only on the coordinates $q_i$. Thus, the 
following relation  has to be valid: $[r_{12}(q),W_1W_2]=0$.
Substituting the explicit form of $r_{12}(q)$ we have
\begin{eqnarray}
\label{rWW}
[r_{12}(q),W_1W_2]_{kl~mn}&=&
W_{kl}\left(\frac{1}{q_{km}}W_{mn}-\frac{1}{q_{km}}W_{kn}-
\frac{1}{q_{ln}}W_{mn}\right)\\
\nonumber
&-&W_{ml}\left(\frac{1}{q_{km}}W_{mn}-
\frac{1}{q_{km}}W_{kn}-\frac{1}{q_{ln}}W_{kn}\right)\\
\nonumber
&+&\delta_{ln}\sum_{j\neq l}\frac{1}{q_{lj}}(W_{kl}W_{mj}-W_{kj}W_{ml}).
\end{eqnarray}

First we show that when $l\neq n$ the first line in (\ref{rWW}) cancels 
the second one.  Since $W_{ij}=\frac{\gamma}{\gamma+q_{ij}}b_j$ we 
get 
\begin{eqnarray} 
\nonumber &&\frac{\gamma}{\gamma+q_{kl}}
\left(\frac{1}{q_{km}}\frac{\gamma}{\gamma+q_{mn}}-\frac{1}{q_{km}}
\frac{\gamma}{\gamma+q_{kn}}-
\frac{1}{q_{ln}}\frac{\gamma}{\gamma+q_{mn}}\right)-\\
\nonumber
&&
\frac{\gamma}{\gamma+q_{ml}}
\left(\frac{1}{q_{km}}\frac{\gamma}{\gamma+q_{mn}}
-\frac{1}{q_{km}}\frac{\gamma}{\gamma+q_{kn}}
-\frac{1}{q_{ln}}\frac{\gamma}{\gamma+q_{kn}}\right)=0.
\end{eqnarray}

In the case $l=n$ the r.h.s. of (\ref{rWW}) reduces to
\begin{eqnarray}
\nonumber
&&-\frac{1}{q_{km}}(W_{kl}-W_{ml})^2+
\sum_{j\neq l}\frac{1}{q_{lj}}(W_{kl}W_{mj}-W_{kj}W_{ml}) =\\
\nonumber
&&-\frac{\gamma^2 q_{km}}{(\gamma+q_{km})^2(\gamma+q_{ml})^2}b_l^2+
\sum_{j\neq l}\frac{1}{q_{lj}}
(\frac{\gamma}{\gamma+q_{kl}}
\frac{\gamma}{\gamma+q_{mj}}
-
\frac{\gamma}{\gamma+q_{kj}}
\frac{\gamma}{\gamma+q_{ml}})b_{j}b_{l} \\
\nonumber
&&
=\frac{q_{mk}}{\gamma+q_{ml}}W_{kl}
\sum_{j}\frac{W_{mj}}{\gamma+q_{kj}}.
\end{eqnarray}
Thus, one has to show that the series
$
S=\sum_{j}\frac{W_{mj}}{\gamma+q_{kj}}
$
vanishes. To this end we consider the following integral
\footnote{We are gratefull to N.A.Slavnov for explaining us
the technique of treating such series.}
$$
I=\frac{1}{2\pi i}\oint \frac{dz}{q_k-z+\gamma}
\frac{\prod_{a\neq m}(q_{a}-z+\gamma)}{\prod_{a}(q_{a}-z)}, 
$$
where the integration contour is taken around infinity. Since the 
integrand is nonsingular at $z\to \infty$, we get $I=0$. On the other
hand, summing up the residues one finds
$
I=\sum_j\frac{1}{\gamma+q_{kj}}
\frac{\prod_{a\neq m} (q_{aj}+\gamma)}{\prod_{a\neq j}q_{aj}}=S.
$

Now we turn to eq.(\ref{qWP}). Explicitly it reads as
\begin{eqnarray}
\nonumber
P_j^{-1}[W_{kl},P_j]&=& 
\left(
\frac{\hbar}{q_{lj}-\hbar}-\frac{\hbar}{q_{kj}}-
\frac{\hbar^2}{q_{kj}(q_{lj}-\hbar)}\right) W_{kl}
+\left(\frac{\hbar^2}{q_{kj}(q_{lj}-\hbar)}+\frac{\hbar}{q_{kj}}\right)W_{jl}
\\
&+& \delta_{jl} \sum_{i\neq j}\left(
\left(
\frac{\hbar^2}{q_{kj}(q_{ij}-\hbar)}-\frac{\hbar}{q_{ij}-\hbar}
\right)W_{ki} -
\frac{\hbar^2}{q_{kj}(q_{ij}-\hbar)}W_{li}\right).
\label{mWP}
\end{eqnarray} 
For the sake of shortness in (\ref{mWP}) we adopt a convention 
that if in some denominator $q_{ij}$ becomes zero, the corresponding 
fraction is also regarded as zero.  Thus, eq.(\ref{mWP}) is
equivalent to the following system of equations
\begin{eqnarray}
\label{knlnj}
P_j^{-1}[W_{kl},P_j] &=& 
\frac{\hbar}{q_{kj}(q_{lj}-\hbar)}(q_{kl}W_{kl}-q_{jl}W_{jl}),
~~\mbox{for}~~k\neq l\neq j;\\
\label{kjnl}
P_j^{-1}[W_{jl},P_j] &=& \frac{\hbar}{(q_{lj}-\hbar)}W_{jl},
~~\mbox{for}~~j\neq l;\\
\label{kln}
P_k^{-1}[W_{kk},P_k] &=& 
-\sum_{i\neq k} \frac{\hbar}{q_{ik}-\hbar}W_{ki}; \\
\label{kjnk}
P_j^{-1}[W_{kj},P_j] &=& 
\frac{\hbar}{q_{kj}}(W_{jj}-W_{kj})+
\frac{\hbar(\hbar-q_{kj})}{q_{kj}}
\sum_{i\neq j}\frac{1}{q_{ij}-\hbar}W_{ki}\\
\nonumber
&-&\frac{\hbar^2}{q_{kj}}\sum_{i\neq j}\frac{1}{q_{ij}-\hbar}W_{ji},~~ 
\mbox{for}~~k\neq j.
\end{eqnarray} 
In the sequel we shall give an explicit proof only for the latter case
since the other three cases are treated quite analogously. The
l.h.s. of (\ref{kjnk}) is
\begin{displaymath}
\label{lh}
P_j^{-1}[W_{kj},P_j] = P_j^{-1}
\frac{\prod_{a\neq k}(q_{aj}+\gamma)}{\prod_{a\neq j}q_{aj}}P_j-W_{kj}=
\gamma
\frac{\prod_{\stackrel{a\neq k}{a\neq j}}(q_{aj}+\gamma-\hbar)}
{\prod_{a\neq j}(q_{aj}-\hbar)}-W_{kj}. 
\end{displaymath}  
As to the r.h.s., one needs to calculate the sum
$\sum_{i\neq j}\frac{1}{q_{ij}-\hbar}W_{ki}$.
For this purpose we evaluate the following integral with the 
integration contour around infinity:
$$
I=\frac{1}{2\pi i}\oint \frac{dz}{z-q_j-\hbar}
\frac{\prod_{a\neq k}(q_{a}-z+\gamma)}{\prod_{a}(q_{a}-z)}, 
$$
The regularity of the integrand at $z\to \infty$ gives $I=0$.  
On the other hand, summing up the residues one finds
$$
I=-
\frac{1}{\hbar}\frac{\prod_{a\neq k}(q_{aj}+\gamma-\hbar)}
{\prod_{a\neq j}(q_{aj}-\hbar)}-
\sum_{i\neq j}\frac{1}{q_{ij}-\hbar}W_{ki}+\frac{1}{\hbar} W_{kj}
$$ 
from here one deduces the desired series:
\begin{equation}
\label{sum}
\sum_{i\neq j}\frac{1}{q_{ij}-\hbar}W_{ki}=
-\frac{1}{\hbar}\frac{\prod_{a\neq k}(q_{aj}+\gamma-\hbar)}
{\prod_{a\neq j}(q_{aj}-\hbar)}+\frac{1}{\hbar} W_{kj}
\end{equation}
and 
\begin{equation}
\label{sum1}
\sum_{i\neq j}\frac{1}{q_{ij}-\hbar}W_{ji}=
-\frac{1}{\hbar}\frac{\prod_{a\neq j}(q_{aj}+\gamma-\hbar)}
{\prod_{a\neq j}(q_{aj}-\hbar)}+\frac{1}{\hbar} W_{jj}.
\end{equation}
Now substituting these sums in the r.h.s. of (\ref{kjnk}) one proves 
(\ref{kjnk}).

It follows from our proof that the $L$-operator
$$
L=\sum_{ij}\frac{\prod_{a\neq i}(q_{aj}+\gamma)}{\prod_{a\neq j}q_{aj}}
e^{-\hbar \frac{\partial}{\partial q_j}}E_{ij}
$$
realizes the representation of the algebra (\ref{qLL}). 

Let us briefly discuss the degeneration of the RS system to the rational CM
model. To get the CM model one should rescale $\hbar \to \gamma\hbar$
and to consider the limit $\gamma\to 0$, $L\to 1+\gamma {\cal L}$.
Then ${\cal L}$ is the $L$-operator of the CM model. From 
eqs.(\ref{QTL}), (\ref{qLL}) one derives the quantum algebra satisfied by 
the $L$-operator of the CM model:  
\begin{eqnarray} 
\label{cmQL} 
[Q_1,{\cal L}_2]&=&\hbar \sum_{i} E_{ii} \otimes E_{ii}, \\
\label{cmLL}
[{\cal L}_1 ,{\cal L}_2]&=&\hbar [r_{12}-\bar{r}_{12},{\cal L}_1]-
\hbar [r_{21}-\bar{r}_{21},{\cal L}_2]+\hbar^2
[r_{12}-\bar{r}_{12},r_{21}-\bar{r}_{21}].
\end{eqnarray}
The last formula can be written in the following elegant form
\begin{equation}
\label{eleg}
[{\cal L}_1+\hbar(r_{21}-\bar{r}_{21}),
{\cal L}_2+\hbar(r_{12}-\bar{r}_{12})]=0.
\end{equation} 

\section{Quantum R-matrix for the trigonometric CM system}
In this section we describe the dual parameterization of $T^{*}G$,
which is related to the trigonometric CM system.
We start with diagonalizing the group element $g=VDV^{-1}$ 
and impose the constraint $Ve=e$. The derivatives 
$\frac{\delta V_{ij}}{\delta g_{kl}}$ and 
$\frac{\delta D_{i}}{\delta g_{kl}}$ are obtained in the same manner
as in Section 2. Calculating the Poisson brackets of $V$, $D$ and $A$
we get
\begin{eqnarray} 
\{V_1, A_2\}&=&V_1V_2s_{12}V_2^{-1}, \label{VA}\\
\{D_1, A_2\}&=&-D_1V_2\sum_i E_{ii}\otimes E_{ii}V_2^{-1}, 
\label{DA} 
\end{eqnarray} 
where 
$$
s_{12}=-\sum_{i\neq j}\frac{D_i}{D_i-D_j}F_{ij}\otimes E_{ji}.
$$

The $\cal L$-operator corresponding to
the trigonometric CM system will be given as the
following function on the phase space:
\begin{equation}
{\cal L}=V^{-1}AV.
\label{Lop}
\end{equation}
Calculation of the Poisson algebra of $T^{*}G$ in terms of 
${\cal L},V,D$ results in 
\begin{eqnarray}                                            
\{{\cal L}_1,{\cal L}_2\}&=&[\tilde{r}_{12},{\cal L}_1]-
[\tilde{r}_{21},{\cal L}_2] 
\label{cLL}\\
\{V_1,{\cal L}_2\}&=&V_1s_{12}, \label{LV}\\
\{D_1,{\cal L}_2\}&=&-D_1\sum_i E_{ii}\otimes E_{ii}, \label{DL}
\end{eqnarray} 
where we have introduced the matrix 
$$
\tilde{r}_{12}=-s_{12}+\frac{1}{2}C.
$$
The matrix $\tilde{r}_{12}$ can be written in the following form
\begin{equation} 
\tilde{r}_{12}=
-\frac{1}{2}\sum_{i\neq j}\mbox{cth} 
\frac{q_{ij}}{2}E_{ij}\otimes E_{ji}+ \frac{1}{2}\sum_{i\neq 
j}\frac{e^{\frac{q_{ij}}{2}}} {\mbox{sinh} \frac{q_{ij}}{2}}E_{ii}\otimes 
E_{ji}+\frac{1}{2}\sum_{i}E_{ii}\otimes E_{ii},  
\end{equation} 
where $q_i=\log D_{i}$.

Now we clarify the connection of $\tilde{r}_{12}$ with the dynamical
$r$-matrix first found in \cite{AT}. By conjugating $\cal L$-operator
(\ref{Lop}) with the matrix $D^{\frac{1}{2}}$: 
$\tilde{\cal L}=D^{\frac{1}{2}}{\cal L}D^{-\frac{1}{2}}$ and 
calculating the Poisson bracket for $\tilde{\cal L}$ with the help of 
(\ref{cLL}) and (\ref{DL}) we arrive at 
\begin{equation} 
\{\tilde{\cal L}_1,\tilde{\cal L}_2\}=[\tilde{R}_{12},\tilde{\cal L}_1]- 
[\tilde{R}_{21},\tilde{\cal L}_2], 
\label{LLL} 
\end{equation} 
where 
\begin{eqnarray}
\nonumber
\tilde{R}_{12}&=&D^{\frac{1}{2}}_{1}D^{\frac{1}{2}}_{2}\tilde{r}_{12}
D^{-\frac{1}{2}}_{1}D^{-\frac{1}{2}}_{2}
-\frac{1}{2}\sum_{i}E_{ii}\otimes E_{ii}\\
&=&
-\frac{1}{2}\sum_{i\neq j}\mbox{cth} 
\frac{q_{ij}}{2}E_{ij}\otimes E_{ji}+ 
\frac{1}{2}\sum_{i\neq j}\frac{1}{\mbox{sinh} \frac{q_{ij}}{2}}E_{ii}\otimes E_{ji}.  
\label{new}
\end{eqnarray}
It is important to note that  $\tilde{R}$ differs from
the matrix found in \cite{AT} by the term 
$$
\frac{1}{2}\sum_{i\neq j}
\frac{1}{\mbox{sinh} \frac{q_{ij}}{2}}E_{ii}\otimes (E_{ji}+E_{ij}).  
$$
However, this term does not contribute to the bracket (\ref{LLL})
if we take into account the representation of the 
$\tilde{\cal L}$-operator of the trigonometric CM system:
\begin{equation}
\tilde{\cal L}=\sum_{i}p_iE_{ii}+\frac{1}{2}\sum_{i\neq 
j}\frac{1}{\mbox{sinh}\frac{q_{ij}}{2}}E_{ij},
\label{cmL}
\end{equation}   
where $(p,q)$ form a pair of canonically conjugated variables.
Thus, on the reduced space these matrices define the same Poisson 
structure for $\tilde{\cal L}$.

Now we quantize $T^*G$ in terms of $A$, $V$ and $D$ variables.
We postulate the following commutation relations
\begin{eqnarray} 
\left[ V_1, A_2\right] &=&\hbar V_1V_2s_{12}V_2^{-1}, 
\label{qVA}\\ 
\left[ D_1, A_2\right] &=&-\hbar D_1V_2\sum_i 
E_{ii}\otimes E_{ii}V_2^{-1}.  
\label{qDA} 
\end{eqnarray} 
One can check that the compatibility of these relations with (\ref{qAg}) 
follows from the following identity satisfied by $s_{12}$:  $$ 
s_{12}-D_1^{-1}s_{12}D_1+\sum_{i}E_{ii}\otimes E_{ii}=C.
$$
 
By using eqs.(\ref{qAA}), (\ref{qVA}) and (\ref{qDA}) one derives
the commutation relations for the quantum $\cal L$-operator:
\begin{eqnarray}                                            
\left[ {\cal L}_1,{\cal L}_2\right]&=&\hbar [\tilde{r}_{12},{\cal L}_1]-
\hbar[\tilde{r}_{21},{\cal L}_2]+\hbar^2 
[ \tilde{r}_{12},\tilde{r}_{21}] 
\label{qcmLL}\\ 
\left[ D_1,{\cal L}_2\right] &=&-\hbar D_1\sum_i E_{ii}\otimes E_{ii},  
\label{qDL} \\
\label{qVL}
\left[V_1,{\cal L}_2\right]&=&\hbar V_1 s_{12}.
\end{eqnarray} 
The relation (\ref{qcmLL}) can be also written in the form (\ref{eleg}): 
\begin{equation}
[{\cal L}_1+\hbar \tilde{r}_{21},
{\cal L}_2+\hbar \tilde{r}_{12}]=0.
\label{sg}
\end{equation} 
One can check without problems that the $L$-operator (\ref{cmL})
realizes the representation of the algebra (\ref{sg}).

To complete our discussion let us show the existence of $N$ mutually
commuting operators in the algebra composed by the $L$-operator and
the coordinates $D$. Obviously, $I_n=\mbox{tr}A^n$ mutually
commute. Applying the technique used in the previous section to derive
the quantum integrals of motion, one can show
that $I_n$ can be expressed as the following function 
of $\cal L$ and $D$:
$$
I_n=\mbox{tr}A^n=\mbox{tr}V{\cal L}^nV^{-1}=
\mbox{tr}_{12}C_{12}^{t_2}({\cal L}_1+\hbar s_{21}^{t_2})^n. 
$$
In the component form these integrals look as
$$
I_n=
\left({\cal L}_{j_1j_2}\delta_{j_1m_1}+ 
\hbar\frac{\delta_{j_1j_2}-1}
{D_{j_1}D_{j_{m_1}}^{-1}-1}\right)\ldots
\left({\cal L}_{j_nj_{n+1}}\delta_{j_nm_n}+ 
\hbar\frac{\delta_{m_{n-1}j_{n+1}}-\delta_{m_{n-1}j_n}}
{D_{j_n}D_{j_{m_n}}^{-1}-1}\right)\delta_{j_{n+1}m_{n}}.
$$
Let us note that $I_n$ can not be expressed as a linear combination
of $\mbox{tr}{\cal L}^n$ solely.

\section{Conclusion}
The approach for $R$-matrix quantization of the RS models proposed in 
the paper seems to be general. The problem of a real interest is to 
apply it to the trigonometric and elliptic cases. 
As was recently shown the trigonometric and 
elliptic RS models are obtained from the cotangent bundles over the 
centrally extended loop group \cite{GN} and double loop group 
\cite{AFM1} respectively. The natural suggestion is to use for this
purpose the above-mentioned phase spaces.

It is known that the Heisenberg double \cite{Sem} can be regarded as 
a natural deformation of the cotangent bundle $T^*G$. It seems to be
interesting to investigate the Poisson structure of the Heisenberg
double in the same parameterization. One could expect the appearance
of another Poisson structure on the Frobenius group induced 
by the one on the dual Poisson-Lie group $G^*$.
 
The appearance of the quantum Frobenius group $F$ states the problem
of developing the corresponding representation theory. Owing to the
method of orbits one can suggest that irreducible representations of 
$F$ should be in correspondence with the symplectic leaves of the 
Poisson-Lie structure. On the other hand, it is known \cite{SKyoto} that 
the symplectic leaves of a Poisson-Lie structure are the orbits of the 
dressing transformation. Studying the orbits of $F^*$
and corresponding representations of $F$ one can hope to obtain the quantum 
integrable systems being some "spin" generalizations of the RS model 
\cite{KrZ}. Another open problem related to the representation
theory is to find the universal Frobenius $R$-matrix.

As was shown in \cite{H} the Ruijsenaars Hamiltonians can be 
related to the special $L$-operator satisfying the fundamental
relation $RLL=LLR$ with Belavin's elliptic $R$-matrix. It would
be interesting to clarify the relationship of this approach 
with our construction.

{\bf ACKNOWLEDGMENT} The authors are 
grateful to L.O.Chekhov, P.B.Medvedev and N.A.Slavnov for valuable 
discussions.  This work is supported in part by the RFFR grants 
N96-01-00608 and N96-01-00551 and by the ISF grant a96-1516.  
 

\begin{thebibliography}{99}
{\small
\bibitem{AT} J.Avan and M.Talon, Phys.Lett.B303 (1993) 33-37.
\bibitem{BV} O.Babelon and C.M.Viallet, Phys.Lett.B237 (1989) 411.
\bibitem{ABT} J.Avan, O.Babelon and M.Talon, Alg.Anal. 6(2) (1994) 67.
\bibitem{Skl} E.K.Sklyanin, Alg.Anal.,  6(2) (1994) 227.
\bibitem{BS} H.W.Braden and T.Suzuki, Lett.Math.Phys. 30 (1994) 147.
\bibitem{AR} J.Avan and G.Rollet, The classical r-matrix for the 
relativistic Ruijsenaars-Schneider system, preprint BROWN-HET-1014 (1995).  
\bibitem{Sur} Yu.B.Suris, Why are the rational and hyperbolic
Ruijsenaars-Schneider hierarchies governed by the same $R$-matrix
as the Calogero-Moser ones ? hep-th/9602160.
\bibitem{NKSR} F.W.Nijhoff, V.B.Kuznetsov, E.K.Sklyanin and 
O.Ragnisco, Dynamical $r$-matrix for the elliptic Ruijsenaars-Schneider 
model, solv-int/9603006.  
\bibitem{Surel} Yu.B.Suris, Elliptic Ruijsenaars-Schneider and 
Calogero-Moser hierarchies are governed by the same $r$-matrix,
solv-int/9603011.
\bibitem{R} S.N.Ruijsenaars, Comm.Math.Phys. 110 (1987) 191.
\bibitem{ABB} J.Avan, O.Babelon and E.Billey, The Gervais-Neveu-Felder
equation and the quantum Calogero-Moser systems, preprint PAR LPTHE 95-25,
May 1995; hep-th/9505091 (to appear in Comm.Math.Phys.).  
\bibitem{GN} J.L.Gervais and A.Neveu, Nucl.Phys. B238 (1984) 125.
\bibitem{Fel} G.Felder, Conformal field theory and integrable systems
associated to elliptic curves, hep-th/9407154.
\bibitem{BBB} O.Babelon, D.Bernard and E.Billey, A quasi-Hopf algebra 
interpretation of quantum 3-j and 6-j symbols and difference equations,
preprint PAR LPTHE 95-51, IHES/P/95/91, q-alg/9511019.
\bibitem{AM} G.E.Arutyunov and P.B.Medvedev, 
Generating equation for $r$-matrices related to dynamical
systems of Calogero type, hep-th/9511070, to appear in Phys.Lett.A.
\bibitem{AFM} G.E.Arutyunov, S.A.Frolov and P.B.Medvedev,
Elliptic Ruijsenaars-Schneider model via the Poisson reduction
of the affine Heisenberg double, hep-th/9607170.
\bibitem{AFM1} G.E.Arutyunov, S.A.Frolov and P.B.Medvedev
Elliptic Ruijsenaars-Schneider model from the cotangent bundle
over the two-dimensional current group, hep-th/9608013.
\bibitem{AF} A.Alekseev and L.D.Faddeev, 
Comm.Math.Phys. 141 (1991) 413.
\bibitem{GN} Gorsky A. and Nekrasov N.,  Nucl.Phys. B414 (1994) 213;
Nucl.Phys. B436 (1995) 582; A.Gorsky, 
Integrable many body systems in the field theories, 
Prep. UUITP-16/94, (1994).
\bibitem{F} L.D.Faddeev Integrable models in (1+1)-dimensional
quantum field theory. - In "Recent advances in field theory and 
statistical mechanics". Eds. Zuber J.B., Stord R. (Les Houches Summer
School Proc. session XXXiX, 1982), Elsvier Sci.Publ., (1984) 561.
\bibitem{KS} P.P.Kulish, E.K.Sklyanin, Quantum spectral transform 
method. Recent developments. - In "Integrable quantum field 
theories". Eds. Hietarinta J., Montonen C., Lect.Not.Phys. 51 
(1982) 61.
\bibitem{OP1} M.A.Olshanetsky, A.M.Perelomov, Phys. Reps. 71 (1981) 313.  
\bibitem{BD} A.A.Belavin and V.G.Drinfel'd, 
Funk.Anal.i ego pril. 16(3) (1982) 1-29.
\bibitem{Chari} V.Chari and A.Pressley, A Guide to Quantum Groups,
Cambridge University Press.
\bibitem{M}  J.M.Maillet, Phys.Lett.B245 (1990) 480.  
\bibitem{Sem} M.A.Semenov-Tian-Shansky, Teor.Math.Phys. 93 
(1992) 302 (in Russian).  
\bibitem{SKyoto}  M.A.Semenov-Tian-Shansky, 
Publ.RIMS Kyoto Univ. 21(6) (1985) 1237-1260.  
\bibitem{KrZ} I.Krichever, A.Zabrodin, Spin 
generalizations of the Ruijsenaars-Schneider model, non-abelian 2D Toda 
chain and representations of Sklyanin algebra, hep-th/9505039.  
\bibitem{H} K.Hasegawa, Ruijsenaars' commuting difference operators
as commuting transfer matrices, q-alg/9512029.
} 
\end{thebibliography}
\end{document}